\documentclass[12pt]{amsart}

\usepackage{amsfonts}
\usepackage{amsmath}
\usepackage{amssymb}
\usepackage{amscd, amsthm}
\usepackage{verbatim}
\usepackage{latexsym}

\newcommand{\N}{\mathbb{N}}
\newcommand{\T}{\mathbb{T}}
\newcommand{\A}{\mathbb{A}}
\newcommand{\Q}{\mathbb{Q}}

\newcommand{\C}{\mathbb{C}}

\newcommand{\aQ}{\overline{\mathbb{Q}}}
\newcommand{\Z}{\mathbb{Z}}

\newcommand{\aFp}{\overline{\mathbb{F}_{p}}}

\newcommand{\aQp}{\overline{\mathbb{Q}_{p}}}
\newcommand{\aQl}{\overline{\mathbb{Q}_{\ell}}}

\newcommand{\GL}{\mathrm{GL}}

\newcommand{\rhobar}{\overline{\rho}}
\newtheorem{theorem}{Theorem}[section]
\newtheorem{lemma}[theorem]{Lemma}
\newtheorem{prop}[theorem]{Proposition}

\newtheorem{cor}[theorem]{Corollary}

\newcommand{\Ad}{{\rm Ad}}

\def\rhobar{ {\bar {\rho} } }
\def\rhob{ {\bar {\rho} } }

\newcommand{\Galois}{\mathrm{Gal}}
\newcommand{\Gal}{\Galois( \bar{ {\mathbb Q}}/{\mathbb Q})}
\newcommand{\F}{{\mathbb F}}

\newcommand{\vep}{{\varepsilon}}

\begin{document}

\title[The level 1 case of Serre's conjecture]
{On Serre's modularity conjecture for 2-dimensional mod $p$
 representations of $\Gal$ unramified outside $p$}
\author[Chandrashekhar Khare]{Chandrashekhar Khare}
\email{shekhar@math.utah.edu}
\address{Department of Mathematics \\
         155 South 1400 East, Room 233 \\
           Salt Lake City, UT 84112-0090 \\
         U.S.A.}


\maketitle

\centerline{\it To my father on completing 50 years of excellence}

\begin{abstract}
We prove the level one case of Serre's conjecture. Namely, we prove that
any
continuous, odd, irreducible representation $\rhobar:\Gal \rightarrow
GL_2(\aFp)$ which is unramified outside $p$ arises from a cuspidal eigenform in
$S_{k(\rhobar)}(SL_2(\Z))$. The proof relies on the methods introduced in 
an earlier joint work with J-P.~Wintenberger, together with a new method 
of ``weight reduction''. 
\end{abstract}

\tableofcontents

\section{Introduction}

Fix a continuous, absolutely irreducible, 2-dimensional, odd, mod $p$ representation
$\rhob:\Gal \rightarrow \GL_2(\F)$ with $\F$ a finite field of
characteristic $p$.
We denote by $N(\rhobar)$ the (prime to $p$) Artin conductor of $\rhobar$, and 
$k(\rhobar)$ the weight of $\rhobar$ as defined in \cite{Serre3}. Serre has conjectured in \cite{Serre3} that such a $\rhob$ {\it arises} (with respect to some fixed embedding $\iota:\aQ \hookrightarrow \aQp$) from a newform of weight $k(\rhobar)$ and level $N(\rhobar)$. 

If $\rhob$ is unramified outside $p$, we say that it is of level 1.
This corresponds to the case 
when $N(\rhobar)$ is $1$.

\subsection{The main result}

\begin{theorem}\label{main}
  A $\rhob$ of level one arises from $S_{k(\rhob)}(SL_2(\Z))$ with respect
  to an embedding $\iota:\aQ \hookrightarrow \aQp$.
\end{theorem}

The theorem settles the conjecture stated 
in article 104 of \cite{Serre2} which is in the case of level 1.
We summarise the history of the conjecture. The level 1  
conjecture was first 
made by Serre in September 1972 (just after the Antwerp meeting)  when he
wrote to Swinnerton-Dyer about it, and asked him to mention this
conjecture (or problem) in his Antwerp text (see \cite{SWD}, p.9). This is
probably the first appearance (1973) of this conjecture in print.
Serre wrote about these conjectures to Tate on May 1st, 1973.
Tate replied to Serre first on June 11, and then on July 2, 1973: in the second letter, he 
proved the conjecture for $p = 2$.

The proof of Theorem \ref{main} builds on the ideas of an earlier work of
Wintenberger and the author, see \cite{KW}. There the above
theorem was proved for primes $p=5,7$ (it being known 
earlier conditionally under GRH for $p=5$ by \cite{Brue}), it being already known
for $p=2,3$ because of a method of Tate, see \cite{Tate} for $p=2$, which was 
later applied to the case of $p=3$ by Serre, see page 710 of
\cite{Serre2}. In \cite{KW} it was shown how
modularity lifting theorems would yield Serre's conjecture when proved in
sufficient generality, and the conjecture was proven in level 1 for
weights $2,4,6,8,12,14$.
The main contribution of this paper is a method to 
prove the level 1 case of the
conjecture using only known modularity lifting theorems, thus completing the proof of the level 1 case of Serre's conjecture. We need lifting theorems when either the $p$-adic lift is crystalline at $p$ of weight $k$ (i.e., Hodge-Tate weights $(k-1,0)$) $\leq p+1$ (and when the weight is $p+1$ the lift is ordinary at $p$), or at $p$ the lift is of Hodge-Tate weights $(1,0)$ and 
Barsotti-Tate over $\Q_p(\mu_p)$.

Henceforth $p$ will be an odd prime.
We use the inductive method proposed in Theorem 5.1 of \cite{KW}
to prove the level 1 case of the conjecture. 
The main new idea of this paper is a ``weight reduction'' 
technique. This  allows us to carry out 
the inductive step in a manner that is different from the one contemplated in 
{\it loc.\ cit.\ }

The cases of the conjecture for small weights in level 1 proved in \cite{KW}
were dealt with using the
results of Fontaine, Brumer and Kramer, and Schoof (\cite{Fontaine}, \cite{BK} and \cite{Schoof1}), together with modularity
lifting results of the type pioneered by Wiles.  

In this paper we use the results of \cite{KW} for weights $2,4,6$, and after that prove the level
one case without making any further use of results classifying abelian varities
over $\Q$ with certain good reduction properties. Thus in the end we see that the
only such results we use are those showing that there is no semistable
abelian variety
over $\Q$ with good reduction outside 5 (see \cite{Fontaine}, \cite{Schoof1} and \cite{BK}). 
We do make use of course of modularity lifting results. These are 
due to Wiles, Taylor, Breuil, Conrad, Diamond, Fujiwara, Kisin, Savitt, Skinner et al
(see \cite{Wiles}, \cite{TW},  
\cite{Fujiwara}, \cite{CDT},\cite{BCDT}, \cite{SW2}, \cite{SW3}, \cite{Savitt}, \cite{Kisin}). In particular, besides the basic method 
of Taylor and Wiles in \cite{Wiles} and \cite{TW}, we need crucially the results of Skinner and Wiles in \cite{SW2},
\cite{SW3}, and the result of Kisin in \cite{Kisin}.  Although Kisin proves a very general modularity lifting theorem for potentially Barsotti-Tate lifts (at $p$) when the residual representation is non-degenerate, i.e.,  
irreducible on restriction to $G_{\Q(\sqrt{{ (-1)}^{p-1 \over 2}p})}$, in this paper the main theorem of \cite{Kisin} is used only in the 
case when the $p$-adic lift being considered is (locally at $p$) Barsotti-Tate over $\Q_p(\mu_p)$. The residually degenerate cases are handled by quoting the results of Skinner and Wiles in \cite{SW2}, \cite{SW3} which may be applied as in these cases the lifts that need to be proved modular are ordinary up to a twist. The ordinarity is  a consequence of a  result
of Breuil and M\'ezard (Proposition 6.1.1 of \cite{BM}), and Savitt (see \cite{Savitt1}, Theorems 6.11 and 6.12), which is vital for us. 

As in \cite{KW}, we use crucially the
potential version of Serre's conjecture proved by Taylor in \cite{Tay1} and \cite{Tay2}, and a deformation theoretic result of B\"ockle in the appendix to \cite{[K03]}.

Theorem \ref{main} yields the following corollaries (see Section \ref{cors}),
the first needing also the method of ``killing ramification'' of Section 5.2 of \cite{KW}:

\begin{cor}\label{cond}
  If $\rhob$ is  an irreducible, odd, 2-dimensional, mod $p$ representation of $G_\Q$ with $k(\rhobar)=2$, 
  $N(\rhobar)=q$, with $q$ prime, and $p>2$, then
  it arises from $S_2(\Gamma_1(q))$.
\end{cor}

\begin{cor}\label{finite}
 There are only finitely many isomorphism classes of continuous semisimple
 odd representations $\rhob:G_\Q \rightarrow GL_2(\aFp)$
 that are unramified outside $p$.
\end{cor}

Our theorem, when combined with modularity lifting theorems also implies that if $\rho:G_\Q \rightarrow GL_2({O})$ is an irreducible $p$-adic representation unramified outside $p$ and at $p$ crystalline of Hodge Tate weights $(k-1,0)$ with $k$ even and 
either $2 \leq k \leq p-1$ (even the weight $p+1$ case can be deduced, after the work in \cite{BLZ}: 
see Lemma \ref{rubbish} below) or $\rho$ is ordinary at $p$, then $\rho$ arises from $S_k(SL_2(\Z))$. It is quite likely that the restriction on weights in the non-ordinary cases
can be eased to allowing weights up to $2p$ provided that residually the representation is non-degenerate. This will probably follow
from ongoing work of Berger on Breuil's conjecture in \cite{Breuill}, and the modifications of the Taylor-Wiles system carried out in \cite{Kisin}.
Theorem \ref{main} also implies that
a $\GL_2$-type semistable abelian variety over $\Q$ with good reduction outside a prime $p$ is a factor of $J_0(p)$. Such a result was earlier used in \cite{KW} in the case when $p=5,7,11,13$ being a special case of the results proven in \cite{Schoof1}, \cite{BK}: now we can recover these results in the case of $\GL_2$-type abelian varieties when $p>5$. 

\subsection{Sketch of proof}

We give a 
rough sketch of the proof of Theorem \ref{main}, starting with some general comment about the method used. As in Theorem 5.1 of \cite{KW}, the method is inductive with respect to the prime which is the residue characteristic, but as said earlier the inductive step is carried out differently. The method uses in an essential way the method of 
``congruences between Galois representations'' which was introduced in Section 4 of \cite{KW}
to prove the cases of the conjecture in level 1 and weights $6,8,12,14$, and is refined here. In Section 4 of {\em loc.\ cit.\ }   congruences were produced
between $p$-adic representations of $\Gal$ 
that were crystalline of weight $p+1$ at $p$ and semistable of weight 2 at $p$
(the analog for modular forms being a result of Serre: see Th\'eor\`eme 11 of 
article 97 of \cite{Serre2}).
The  method in Section 2 of \cite{KW} can be used to prove more results about such congruences which parallel 
results that are well-known for congruences between modular forms. For instance
we can now in principle prove analogs for Galois representations of 
the ``type changing'' arguments of Carayol in \cite{car} for modular forms (some instances are carried out in 
Section 3 of the paper, and used in the proof of Theorem \ref{main} to ``change types'' at a prime different from the residue characteristic), or the level raising results for
modular forms of Ribet (see  \cite{Ribet1}). We do not use the Galois-theoretic
analog of the latter in this paper, but this and other such ``level raising'' results for Galois representations will be crucial in future work. 

We fix an embedding $\iota_p:\aQ \hookrightarrow \aQp$ for each prime $p$. 
We say a residual mod $p$ representation $\rhobar$ or a $p$-adic 
representation $\rho$, is {\it modular} if it is either reducible (in the $p$-adic case we assume irreducibility) or it arises from a newform with respect to the embedding $\iota_p$.
We say that a compatible system $(\rho_{\lambda})$ is modular if it arises from a newform.

Assuming we have proved Serre's conjecture for level 1 modulo the $n$th prime $p_n$, we prove it mod $P_{n+1}$ where $P_{n+1}$ is the least non-Fermat prime $>p_n$. This is (more or less!) the inductive method 
to attack the level 1 case of Serre's conjecture 
proposed in Theorem 5.1 of \cite{KW}.
Via the methods of \cite{KW}, see also Lemma \ref{rubbish} and Corollary \ref{trivial} below, this also means that one knows (the level 1 case of) Serre's conjectures for $\rhobar$ of any residue characterictic $p$ bigger than $p_n$ for all weights up to $P_{n+1}+1$. Then we repeat the process starting with $P_{n+1}$ instead of $p_n$.

In \cite{KW} to prove weights $6,8,12,14$ the minimal lifting result in Section 2 of \cite{KW} was used twice to get 2 different compatible systems whose interplay (via the residual representations at a certain place of the 2 compatible systems being isomorphic up to semisimplification: we say that 2 such compatible systems are {\em linked}) proved the modularity
of $\rhobar$. Here 3 compatible systems are considered instead which arise from $\rhobar$ directly or indirectly.
The lifts constructed are not always minimal. The inductive step is different
from the one proposed in \cite{KW} 
in that the residual modularity is essentially used only for representations in  smaller weights of the same residual characteristic which has already been 
inductively established earlier. Another prime $\ell$ is used as a foil in an auxiliary fashion to achieve this reduction of weight.

More precisely, starting with a $\rhobar$ mod $p:=P_{n+1}$ such that
$p_{n}+1 < k(\rhobar)\leq p+1$ of level 1, we first reduce to considering 
$\rhobar$ which are ordinary at $p$ (by Lemma \ref{weights} below), and then construct 
a minimal lifting (see Proposition \ref{p} below) of $\rhobar$ to a $p$-adic representation that is unramified outside $p$ and of weight 2 at $p$ (it is 
Barsotti-Tate 
over $\Q_p(\mu_p)$ if $k(\rhobar) \neq p+1$, and otherwise semistable of weight 2). We then get a compatible system $(\rho_{\lambda})$ (see Proposition \ref{c}) such that $\rho$ is
part of this system for a $\lambda$ above $p$. 
We choose an odd prime $\ell$ such that $\ell^r||p-1$ (hence the assumption that $p$ is not a Fermat prime), 
and consider the residual representation $\overline \rho_{\lambda}$ for $\lambda$ now a prime above $\ell$. 

Now in the cases when the residual modularity is not known for this mod $\ell$ representation 
(the residual modularity will be known if the image is solvable, or the representation is unramified at $p$) we 
construct another lifting (see Proposition \ref{q} which plays in some sense the analog of a lemma of Carayol \cite{car}, or 
de Shalit's lemma as in \cite{Wiles} and \cite{TW}, for Galois representations) 
$\rho'$ of this mod $\ell$ representation that is Barsotti-Tate at $\ell$, and at $p$ is non-minimal with 
``nebentype'' that is well-chosen, and is unramified outside $\ell,p$. 
(If the mod $\ell$ modularity were known, we conclude the modularity of $(\rho_\lambda)$ 
by a modularity lifting theorem applied at $\ell$: the modularity lifting result in 
this case would be the one contained in \cite{Wiles}, \cite{TW}, \cite{SW2} and \cite{SW3}.)
We construct another compatible 
system $(\rho'_{\lambda})$ of which $\rho'$ is a member, and such that at 
a prime above $p$, because of the well-chosenness of the nebentype at $p$ in the $\ell$-adic lift, the residual representation 
is already known to be modular by the inductive hypothesis (which includes the case when the representation 
is reducible as these by convention are also called modular). This control of the weight of the mod $p$ 
residual representation arising from $(\rho'_{\lambda})$ is due 
to a result of Breuil and M\'ezard in \cite{BM}, and Savitt \cite{Savitt1}. Then the fact that we are 
in a position to apply known modularity lifting results (see \cite{SW2}, \cite{SW3}, \cite{Kisin}) 
is again because of \cite{BM} which says that all lifts of reducible 2-dimensional mod $p$ representation of $G_{\Q_p}$
that become Barsotti-Tate over $\Q_p(\mu_p)$ are up to twist ordinary. (This is one important reason why the 
strategy here does not need modularity lifting results beyond the known range.)

These known modularity lifting theorems then prove modularity of 
$(\rho'_{\lambda})$, and then another use of modularity lifting theorems (of \cite{Wiles}, \cite{TW}, \cite{SW2}, \cite{SW3}) 
prove the modularity of $(\rho_{\lambda})$, as the 2 compatible systems are linked mod the prime above $\ell$ fixed by $\iota_\ell$, 
and thus $\rhobar$ is modular. We need a third kind of compatible system $(\rho''_{\lambda})$
constructed in \cite{KW} to conclude now that all level 1 representations modulo a prime $\geq k(\rhobar)-1$ 
and of the weight of $\rhobar$ are modular (see Lemma \ref{rubbish} and 
Corollary \ref{trivial} below).

The entire strategy of the paper roughly uses that $p_n$ 
is at least two-thirds as large as $p_{n+1}$ (using classical Chebyshev estimates: 
see Section \ref{cheb}). 

Thus schematically the argument may be summarised as follows:
\begin{itemize}

\item Start with a mod $p$ representation
$\rhobar$,

\item lift it to $\rho$ (using Proposition \ref{p}), 

\item construct a compatible system
$(\rho_{\lambda})$ (using Proposition \ref{c}) such that for the prime above $p$ fixed by $\iota_p$ the corresponding representation
is $\rho$, 

\item consider the residual representation at a prime $\lambda$ above 
a ``suitable'' prime $\ell$ ($\lambda$ fixed by $\iota_\ell$),

\item lift this to a  ``good'' lift  $\rho'$  (using Proposition \ref{q}), 

\item make it part of a compatible system 
$(\rho'_{\lambda})$ (using Proposition \ref{p}),  

\item consider the residual representation at the above $p$ fixed by $\iota_p$
arising from the system $(\rho'_{\lambda})$, 

\item inductively this is known to be modular as $\rho'$ is a ``good'' lift, 

\item modularity lifting theorems imply $(\rho'_{\lambda})$ is modular, 

\item another application of lifting theorems gives  $(\rho_{\lambda})$ is modular, and hence $\rhobar$ is modular. 
\end{itemize}

\noindent (This is in the ``generic'' case, as sometimes the procedure yields success earlier.) After this we use a third kind of compatible system 
$(\rho''_{\lambda})$ constructed in Sections 2 and 3 of \cite{KW}, to deduce that once Serre's conjecture is known in weight $k$ for a prime $p$ then it is known in weight $k$ for all primes $\geq k-1$ (see Corollary \ref{trivial} below).

Just as in \cite{KW}, besides modularity lifting results, the potential version of Serre's conjectures 
proven by Taylor in \cite{Tay1} and \cite{Tay2}, and a result of B\"ockle
in the appendix to \cite{[K03]}, is crucial in constructing (minimal and non-minimal) liftings 
$\rho$ of various kinds (see Sections 2 and 3), and then making them part of a compatible system 
$(\rho_\lambda)$ (see Section 4)  whose refined properties are then obtained by arguments 
of Dieulefait and Wintenberger,
\cite{D} and \cite{Wint}.


\section{Liftings}

For $F$ a field, $\Q \subset F \subset \aQ$,
we write $G_{F}$ for the Galois group of $\aQ / F$.
For $\lambda$ a prime/place of $F$,
we mean by $D_{\lambda}$ (resp., $I_{\lambda}$) a decomposition
(resp., inertia) subgroup of $G_F$ at  $\lambda$.
We have fixed embeddings $\iota_p, \iota_\infty$ of $\aQ$ in its completions $\aQp$ and $\C$ in the introduction.
Denote by $\chi_p$ the $p$-adic
cyclotomic character, and $\omega_p$ the Teichm\"uller lift of the mod $p$ cyclotomic character $\overline \chi_p$ 
(the later being the reduction mod $p$ of $\chi_p$). By abuse of notation we also denote by $\omega_p$ the $\ell$-adic character $\iota_\ell\iota_p^{-1}(\omega_p)$ for any prime $\ell$: this should not cause confusion as from the context it will be clear where the character is valued. For a number field $F$ we denote the restriction of a character of $\Gal$  to $G_F$ by the same symbol.

Let $p$ an odd prime.
Fix  $\rhob: G_{\Q}\rightarrow \mathrm{GL}_2 (\aFp )$
to be an odd irreducible representation.
We assume that the Serre
weight $k(\rhob)$ is such that $2 \leq k(\rhob) \leq p+1$. (Note that
there
is always
a twist of $\rhob$ by some power of the mod $p$ cyclotomic character
$\overline{\chi _p}$
that has weights in this range.) We denote by $\Ad^0(\rhobar)$ the $G_{\Q}$-module arising from the adjoint action on the trace 0 matrices of $M_2(\F)$.

Let $\F \subset \aFp$ be a finite field
such that the image of $\overline{\rho}$ is contained in
$\mathrm{GL}_2 (\F )$, and let $W$ be the Witt vectors $W(\F )$.
By a {\it lift} of $\overline{\rho}$, we mean a continuous representation
$\rho : G_{\Q}\rightarrow
\mathrm{GL}_2 ({V})$, where $ V$ is the ring of integers of a
finite
extension of the field of fractions of  $W$, such that
the reduction of $\rho$ modulo the maximal ideal of ${V}$ is
isomorphic to $\overline{\rho}$. Liftings when considered up to equivalence, i.e., up to conjugation by matrices that are residually the identity, are referred to as deformations.
We say that $\rho$ is minimal at a prime $\ell \neq p$ if it  is minimally
ramified at  $\ell$ in the terminology of \cite{[D97]}.

Let $\vep$ be the Teichm\"uller lift of the  character 
${\rm det}(\rhob)\overline \chi_p^{1-k(\rhob)}$ whose restriction to any open subgroup of $\Gal$ we denote by the same symbol.

\subsection{The method of producing liftings of \cite{KW}}\label{liftings}

In this section we will produce liftings with certain prescribed local properties
of $\rhobar$ (one of the properties being unramified almost everywhere) using the methods of Section 2 of \cite{KW}. We assume 
that $\rhobar$ has non-solvable image and $k(\rhobar) \neq p$.
By Lemma 2.6 of \cite{KW}, this also means that for any totally real field $F$,
 $\rhobar|_{G_F}$ has non-solvable image, and as $2 \leq k(\rhobar) \leq p+1$ and $\neq p$,  if $F$ is unramified at $p$, 
$\wp$ a place of $F$ above $p$, $\rhobar|_{I_{\wp}}$ is non-scalar.

To orient the reader we say a few words about 
the way the method of {\it loc.\ cit.\ } gets used here. We wish to point out that the method there is flexible enough to produce liftings with the desired local properties, provided the local calculations work out.

We 
produce liftings
with  prescribed properties as in Propositions \ref{p} and \ref{q} below, 
by proving, as in \cite{KW},
that a deformation ring $R$ (over a suitable ring of integers $O$ of a finite extension of $\Q_p$, with uniformiser $\pi$)  which parametrises (equivalence classes of) liftings of $\rhobar$
with these prescribed properties is flat 
over $O$. (The prescribed properties will be deformation conditions in the sense of \cite{Mazur} and thus the problem will be representable globally by a 
universal ring $R$.) To do this one shows first that $R/(\pi)$ is finite, which
by Lemma 2.4 of {\it loc.\ cit.\ }, as explained below,
is equivalent to showing that $R_F/(\pi)$ is finite, where $R_F$ is
a certain deformation ring for $\rhobar|_{G_F}$ and $F$ is some totally real
field. 

The finiteness of $R_F/(\pi)$ is established by identifying $R_F$ to
a Hecke algebra $\T_F$ which we know is finite over $\Z_p$. 
A suitable $F$ with this property is produced by using Taylor's results
in \cite{Tay1}, \cite{Tay2}: for instance $F$ is unramified at $p$, and in the cases below $R_F$ can be taken to be a minimal deformation ring (although $R$ need not be a minimal deformation ring, see Proposition \ref{q} for instance).
Here the minimality of $R_F$ at the local defining conditions at 
places $\lambda$ not above $p$ does not need futher comment: for $\lambda$ above $p$ the minimality condition is as explained below the same as the minimal condition at $p$ when defining $R$, and is different in the case
of Propositions \ref{p} and \ref{q}.

To see that $R_F$ can be taken to be the the minimal deformation ring in the 
proofs of Propositions \ref{p} and \ref{q} below we indicate the argument. Denote the universal representation $\rho_R$ and $\rho_{R_F}$ corresponding
to the deformation problem that $R$ and $R_F$ represent (with $R_F$ the minimal deformation ring), and $\overline \rho_R$ and $\overline \rho_{R_F}$ it's reduction mod $\pi$. Then using the properties of the lifts prescribed in 
Propositions \ref{p} and \ref{q} below, for $\ell \neq p$, the order of ${\overline \rho_R} (I_\ell)$ is the same as the order of $\rhobar(I_\ell)$. This is seen by using the proof of the claim after Lemma 2.4 of \cite{KW}. 

This gives that
${\overline \rho_R}|_{G_F}$ is a specialisation of $\overline \rho_{R_F}$. 
The finiteness of $R_F$ as a $O$-module yields that the latter has finite image and hence so does the former, which by Lemma 2.4 of \cite{KW} gives that
$R/(\pi)$ is of finite cardinality. 

We give some more details about the choice of $F$ and the identification of $R_F$ to
a Hecke algebra $\T_F$ (see also proof of Theorem 2.2 of \cite{KW}).
We choose the totally real field $F$, Galois over $\Q$ and of even degree, so that when $\rhobar$ is ordinary at $p$ (resp., supersingular by which we mean locally irreducible at $p$), $F$ is unramified at $p$ 
(resp., split at $p$),
$\rhobar|_{G_F}$ is unramified outside places above $p$, {\it and} such that there is a cuspidal automorphic representation $\pi$ for $GL_2(\A_F)$
that is discrete series of parallel weight $(2,\cdots,2)$ at infinity (resp.,
of parallel weight $(k(\rhobar),\cdots,k(\rhobar))$ at infinity) is unramified at all places not above $p$, and is ordinary at places $\wp$ above $p$ 
of conductor dividing $\wp$ (resp., unramified at places $\wp$ above $p$), and unramified at $p$ when $k(\rhobar)=2$,
that gives rise to $\rhobar|_{G_F}$ with respect to (w.r.t.) the 
embedding $\iota_p$. In the ordinary case (in Proposition \ref{q})
we also use the existence of $F$ unramified at $p$, 
Galois over $\Q$ and of even degree, such that $\rhobar|_{G_F}$ is unramified outside places above $p$, {\it and} such that there is a cuspidal automorphic representation $\pi$ for $GL_2(\A_F)$
that is discrete series of parallel weight $(k(\rhobar),\cdots,k(\rhobar))$ at infinity and is unramified at all finite places (see \cite{Tay1} and Proposition 2.5 of \cite{KW}), that gives rise to $\rhobar|_{G_F}$ with respect to (w.r.t.) the 
embedding $\iota_p$.

The existence of a $F$ with all these properties 
follows from Lemma 1.5 and Corollary 1.7 of \cite{Tay1}, and Proposition 2.5 of \cite{KW}, in the ordinary case, and Theorem 5.7 of \cite{Tay2} in the supersingular case, and uses as an ingredient the level-lowering up to base change method of \cite{SW1}. (The evenness of $[F:\Q]$ is not mentioned in \cite{Tay1} in the ordinary case,
but certainly may be ensured by a further quadratic base change.) 

In the case of Proposition \ref{p}, in which case $F$ 
is unramified at $p$, the minimal deformation ring $R_F$ 
parametrises lifts $\rho$ of $\rhobar|_{G_F}$ that are unramified away from $p$, and
at places $\wp$ above $p$ are such that  $\rho|_{I_{\wp}}$ is of the form 
$$\left( \begin{array}{cc}
    \omega_p^{k-2}\chi_p  & *  \\
    0 & 1
\end{array} \right),$$ if $\rhobar|_{I_{\wp}}$
is of the form $$\left( \begin{array}{cc}
    \overline \chi_p ^{k-1} & *  \\
    0 & 1
\end{array} \right),$$ with $2 \leq k \leq p-1$, with the further condition that if $k(\rhobar)=2$,  $\rho|_{I_{\wp}}$ is Barsotti-Tate, and of determinant (the restiction of) $\vep \omega_p^{k-2}\chi_p$. In the case of Proposition \ref{q}, in which case $F$ may be taken to be split at $p$ when $\rhobar|_{D_p}$ is supersingular, the minimal deformation ring $R_F$ 
parametrises lifts $\rho$ of $\rhobar|_{G_F}$ that are unramified away from $p$, and at places $\wp$ above $p$ the representation is crystalline of weight $k(\rhobar)$, and the lifts have determinant $\vep \eta_q^i \chi_p^{k(\rhobar)-1}$
using the notation of Proposition \ref{q}. (We can also ensure that $\vep|_{G_F} $ or $\vep \eta_q^i|_{G_F}$ is trivial if we want when the Serre weight is even:
the Serre weight will always be even in all applications below.)

It remains to recall the identification of $R_F$ to suitable Hecke algebra $\T_F$.
In the  case of Proposition \ref{p} below, the Hecke algebra $\T_F$ is a ($\Z_p$-)algebra cut out by the 
Hecke action on cusp forms for $GL_2(\A_F)$ that are of weight $(2,\cdots,2)$, 
unramified outside $p$, and at places $\wp$ above $p$ of conductor dividing $\wp$ (and unramified if $k(\rhobar)=2$), and
of central character corresponding to (restriction of) $\vep \omega_p^{k-2}$ by class field theory, with respect to the embedding $\iota_p$.
In the case of Proposition \ref{q} below, the Hecke algebra $\T_F$ is a ($\Z_p$-)algebra 
cut out by the 
Hecke action on cusp forms for $GL_2(\A_F)$ that are of weight $(k(\rhobar),
\cdots,k(\rhobar))$, unramified at all finite places, and of central character 
corresponding to, using it's notation, $\eta_q^i \vep \chi_p^{k(\rhobar)-2}$ by class field theory. 
(The fact that the $\T_F$ is non-zero is a consequence of the results of 
Taylor in \cite{Tay1} and \cite{Tay2} that we have recalled.) 

The identification $R_F\simeq \T_F$
is proved using \cite{Fujiwara} in the ordinary case, and  
Section 3 of \cite{Tay2} in the supersingular case. ({\small 
As \cite{Fujiwara} may not be widely available, note that in the ordinary case,
we are in a situation where we do have a minimal modular lift of $\rhobar|_{G_F}$ that is ordinary at places above $p$, and such that $\rhobar|_{G_F}$ has nonsolvable image, and $F$ is unramified at $p$. Thus the deduction of the isomorphism 
$R_F \simeq \T_F$ is by now standard.})

Note that we are allowing $p=3$, which is a case excluded in
some sections of \cite{Tay2}: thus we say a few words to justify why we still have the results of \cite{Tay2} available. We exploit the fact that we know 
that by our assumption 
$\rhobar|_{G_F}$ is not solvable for any totally real field $F$ (Lemma 2.6 of \cite{KW}), and as
$\overline \chi_3$ restricted to $F$ has order 2, 
 we thus have an auxiliary prime $r$, as guaranteed
by Lemma 3  of \cite{DT} or Lemma 4.11 of \cite{[DDT]}, which handles non-neatness problems. (For instance, this allows us to pass from the results in Section 4 of \cite{Tay2}
to those of Section 5 requiring only that $p>2$, and also to have available the results of Sections 2 and 3 of \cite{Tay2} requiring only that $p>2$.)

The finiteness of $R/(\pi)$ leads to $R$ being  flat (finite, complete intersection) over $ O$ if we know that $R \simeq {O}[[X_1,\cdots,X_r]]/(f_1,\cdots,f_s)$ with $s \leq r$.  By the crucial Proposition 1 of B\"ockle's appendix to \cite{[K03]}, this will follow from some purely local information about the kind of lifts $R$ parametrises. (The local conditions will always be deformation conditions in the sense of \cite{Mazur}.) If $R$ parametrises (equivalence classes of) lifts unramified outside a fixed set of primes, of a certain fixed determinant, only the following local properties at each prime $\ell$ need be checked:  the corresponding local deformation ring $R_{\ell}$ should be a flat, complete intersection  over $O$, of relative dimension 
${\rm dim}{_\F} H^0(D_{\ell},\Ad^0(\rhobar))$ 
when $\ell \neq  p$, and of relative dimension 
${\rm dim}{_\F} H^0(D_{p},\Ad^0(\rhobar))+1$ 
when $\ell =  p$.  We check these local conditions in Propositions \ref{p} and \ref{q} showing that  they follow from results of B\"ockle, Ramakrishna and Taylor, see \cite{Boe}, \cite{Boe1}, \cite{Tay3}, \cite{[R02]} (in all the cases below the local deformation ring turns out to be smooth). At primes $\ell \neq p$ 
where no ramification is allowed the deformation ring is directly checked to be smooth of relative dimension ${\rm dim}{_\F} H^0(D_{\ell},\Ad^0(\rhobar))$. At primes $\ell \neq p$ at which the residual representation is ramified and 
where the corresponding local deformations that are 
allowed are minimal the ring is smooth of the required dimension as checked in Section 3 of \cite{Boe}, \cite{[R02]} (see the 
`` local at $\ell \neq p$'' section) and  \cite{Tay3} (see E1 to E3). 
Thus below we only check the local condition at the residual characteristic 
(where the results are again found in \cite{Boe}, \cite{[R02]} 
and \cite{Tay3}), 
and at a prime $\ell$ where the deformations allowed are not minimal.
(Note that in \cite{Tay3} for the local versal 
deformation rings below, 
at $p$, or when $\ell \neq p$ and the deformations allowed are minimal, 
the unobstructedness of the ring
is checked, and it is shown that the tangent space is of dimension ${\rm dim}{_\F} H^0(D_{\ell},\Ad^0(\rhobar))+\delta_{\ell p}$, which implies that the versal ring is smooth of relative dimension ${\rm dim}{_\F} H^0(D_{\ell},\Ad^0(\rhobar))+\delta_{\ell p}$ over $W$.)

\subsection{Minimal $p$-adic weight 2 
lifts of $\rhobar$}\label{weight2}

We consider (just for this subsection) only
$\rhob$ such that
$\rhob$ is ordinary at $p$, i.e., $\rhob|_{I_p}$ has non-trivial
covariants, and we also assume as before $k(\rhob) \neq p$.
In this subsection we consider lifts 
$\rho$ of $\rhobar$ whose determinant is $\vep\omega_p^{k(\rhob)-2}\chi_p$.

We say that $\rho$ is minimal of weight 2 at $p$ (as in E3 and E4 of
\cite{Tay3})  if 
its determinant is $\vep\omega_p^{k(\rhob)-2}\chi_p|_{D_p}$, and assuming $\rhob|_{I_p}$
is of the form
$$\left( \begin{array}{cc}
    \overline \chi_p ^{k-1} & *  \\
    0 & 1
\end{array} \right),$$ with $2 \leq k \leq p-1$, then $\rho|_{I_p}$ is of
the form
$$\left( \begin{array}{cc}
    \omega_p^{k-2}\chi_p  & *  \\
    0 & 1
\end{array} \right),$$ with the further condition that when $k(\rhobar)=2$, $\rho|_{I_p}$ is Barsotti-Tate. 

If a lift $\rho$ satisfies this condition at $p$ and is minimal at all prime $\neq p$ (and thus necessarily
has determinant fixed as above), then we say that it is minimal
of weight 2.

\begin{prop}\label{p}
Let $p$ be a prime $> 3$.
Let $\overline{\rho} : G_{\Q}\rightarrow
\mathrm{GL}_2 (\F)$ be an odd absolutely irreducible representation.
We suppose that $2 \leq k(\rhob) \leq p+1$ and
 $k(\rhob)\not= p$ and with $\rhob$ ordinary at $p$. Then
$\rhob$ has  a lift $\rho$ that is minimal of weight 2. (Its determinant
is necessarily $\vep\omega_p^{k(\rhob)-2}\chi_p$.)
\end{prop}

\begin{proof} If the image of $\rhob$ is solvable we are done using 
that Serre's conjecture is known in this case even in its refined form (see \cite{Ribet}) and the fact that if $\rhobar$ arises from a mod $p$ ordinary 
eigenform in $S_{k(\rhobar)}(\Gamma_1(N),\aFp)$  ($(N,p)=1$), then it also arises from an ordinary eigenform in
$S_2(\Gamma_1(N) \cap \Gamma_0(p),\aFp({\overline \chi_p}^{k(\rhobar)-2}))$ (see Proposition 8.13 of \cite{Gross} or
Section 6 of \cite{Edix}). So we
now assume that the image of $\rhob$ is not solvable.
  
The proof follows immediately from the method of proof of Theorem 2.1 of
\cite{KW}, as noted in Section \ref{liftings}, on noting the
following local fact: consider the versal deformation ring $R_p$
  that parametrises deformations of $\rhob|_{D_p}$
(to $W$-algebras that are complete Noetherian local (CNL)
rings with residue field $\F$ as usual)
that on inertia $I_p$ have the form $$\left( \begin{array}{cc}
    \omega_p^{k-2}\chi_p  & *  \\
    0 & 1
\end{array} \right),$$ have fixed determinant the image of the character 
$\vep\omega_p^{k-2}\chi_p|_{D_p}$, and in the case of weight $k(\rhobar)=2$ the deformation is Barsotti-Tate. Then $R_p$ is a complete intersection, flat over $W(\F)$
of relative dimension $1+{\rm dim}_{\F}(H^0(D_p,{\rm Ad}^0(\rhobar))$: 
in fact it is even smooth 
of relative dimension $1+{\rm dim}_{\F}(H^0(D_p,{\rm Ad}^0(\rhobar))$ as checked
  in E3 and  E4 of \cite{Tay3} (see also \cite{[R02]}).  This, by Section \ref{liftings}, gives the flatness of the ring $R$ over $W$ that parametrises (equivalence classes of) 
lifts of $\rhob|_{D_p}$ that are minimal of weight 2 and hence we are done.
\end{proof}

\noindent{\bf Remark:} In this paper we use this theorem only for representations unramified outside $p$. We also do not need to use the proposition when the image is solvable. The condition of ordinarity may be removed using
results in \cite{BM}, \cite{Savitt1} (Proposition 6.1.2(iii) of the former, 
Theorem 6.22 of latter), together with the modification of Taylor-Wiles systems
in \cite{Kisin}.

\subsection{A Galois theoretic analog of Carayol's lemma}\label{type}

In this subsection we prove a Galois theoretic analog of Lemme 1 of \cite{car}.
Consider $\rhob: \Gal \rightarrow GL_2(\F)$ that is continuous, odd, irreducible,  $2 \leq k(\rhobar) \leq p+1$ and $k(\rhobar) \neq p$, 
and consider an odd prime $q$ which we assume is ramified in
$\rhob$.
Further assume that $\rhob|_{I_q}$ is of the form
$$\left( \begin{array}{cc}
    \overline \chi  & *  \\
    0 & 1
\end{array} \right),$$ where $\overline \chi$ arises from a mod $p$ character of ${\rm Gal}(\Q_q(\mu_q)/\Q_q)$.
Let $\chi$ be its Teichm\"uller lift. This will be a power of the character $\iota_p\iota_q^{-1}(\omega_q)$
which we recall that by our conventions is again denoted by $\omega_q$.

Assume that $p^r||q-1$ ($r>0$) and consider $\eta_q=\omega_q^{ {q-1} \over p^r }$: this is (for this subsection)
 a character with values in $\aQp^*$.
We denote the corresponding global characters which factor through
${\rm Gal}(\Q(\zeta_q)/\Q)$ by the same symbol. We enlarge $\F$ so that it contains all the ${q-1} \over p^r$th roots of 1.
Below, we denote by $O$ the ring of integers of $W(\F)(\mu_{q-1})$ and by $\pi$ a uniformiser of $O$.

In this section the minimality condition at $p$ we will consider
is of being crystalline of weight $k(\rhobar)$.

\begin{prop}\label{q}
  Let $p$ be an odd prime, fix a $\rhob$ as above (in particular $k(\rhobar)\neq     p$ and $\rhobar|_{I_q}$ has the form above), and assume that $\rhob$ does not have solvable
image. Fix an integer $i$. For some $V$ that is the ring of 
integers of a finite extension of $\Q_p$, there is a $V$-valued lift $\rho$ of $\rhob$ of determinant
$\vep\chi_p^{k(\rhobar)-1}\eta_q^i$, that is minimal at primes outside $p,q$, is minimal
at $p$ (crystalline of weight $k(\rhobar)$), and at $q$, $\rho|_{I_q}$ is of the from
$$\left( \begin{array}{cc}
     \chi \eta_q^i  & *  \\
    0 & 1
\end{array} \right).$$ We say that such a lifting has nebentype $\chi \eta_q^i$ at $q$.
\end{prop}

\begin{proof}
 This follows by the arguments in Section 2 of \cite{KW}, as noted in 
Section \ref{liftings}, 
from the following 2 local
facts: 

- The local deformation ring $R_p$ which parametrises (equivalence classes of)
lifts of $\rhobar|_{G_p}$
to $O$-algebras that are crystalline of weight $k(\rhobar)$, and of fixed determinant, is smooth over $ O$ of dimension ${\rm dim}{_\F} H^0(D_p,\Ad^0(\rhobar))+1$.
This is proved in \cite{[R00]}, \cite{Tay3} (see also discussion in Section 2
of \cite{KW} and Proposition 2.3 of 
\cite{KW} for the case of $k(\rhobar)=p+1$).

- Consider the versal ring $R_q$ that parametrises (equivalence classes of) lifts of
$\rhob|_{D_q}$ to CNL $O$-algebras with residue field $\F$ 
of determinant (the restriction to $D_q$ of) $\vep\chi_p^{k(\rhobar)-1}\eta_q^i|_{D_q}$, and such
that
$\rho|_{I_q}$ is of the from
$$\left( \begin{array}{cc}
     \chi\eta_q^i  & *  \\
    0 & 1
\end{array} \right).$$  That such a ring exists (i.e., the conditions we are defining, which can be interpreted as a condition of ordinarity as we are fixing determinants, are deformation conditions) follows easily from our assumption that $q$ is ramified in $\rhobar$ (see for instance Section 6.2 of \cite{dsl}). The key fact we need is that $R_q$ is a complete intersection, which is flat over $ O$, and of relative
dimension ${\rm dim}_{\F}H^0(D_q,\Ad^0(\rhobar))=1$. (In fact, it is even smooth.) The asserted dimension of the cohomology group is easily checked using
that $\rhobar$ is ramified at $q$.

We now prove that $R_q$ is smooth over $O$ of relative dimension 1.
This essentially follows from Section 2 of \cite{Boe} 
(see Section 2, and in particular Theorem 3.8 and Lemma 3.10: the deformation problem we are describing here 
is one obtained by specializing the $T$ in 
Theorem 3.10 (iii) of {\it loc.\ cit.\ } to a specific value and twisting).
We sketch an argument to be more self-contained. 

The dimension over $\F$  of the mod $\pi$ Zariski tangent space of $R_q$ is 1. This follows from the calculations in Section 1 of \cite{Wiles}, 
see also Section 4.3 of \cite{deshalit} for an exposition in semistable cases, as the dual number lifts that arise from $R_q$ are the same as the minimal lifts in
\cite{Wiles}. Thus to show that $R_q$ is in fact smooth of relative 
dimension 1 over O it will be enough to show that there are infinitely many 
non-equivalent $O$-valued lifts $\rho_q$ of $\rhobar|_{D_q}$ 
of the required kind which we now proceed to show. 
Any lift will be tamely ramified and thus will be specified by lifting the image of $\rhobar(\sigma_q)$ 
and $\rhobar(\tau_q)$, where $\sigma_q,\tau_q$
are generators of Galois group of the maximal tamely ramified extension of $\Q_q$, and the 
only relation these satisfy  is $\sigma_q\tau_q\sigma_q^{-1}=\tau_q^q$. When $\overline \chi$ is non-trivial, 
and thus $\rhobar|_{I_q}$ may be assumed split and 
$\rhobar|_{D_q}$ is diagonal, 
by inspection we get infinitely many lifts to diagonal matrices. 

In the case when $\overline{\chi}$ and $\chi$ are  trivial, and hence $\rhobar|_{I_q}$ is unipotent, and not-trivial by 
assumption, again a simple calculation yields infinitely many $O$-valued lifts.
There will be 2 cases corresponding to $\chi':=\eta_q^i$ being trivial or non-trivial.
When it is trivial this is covered by E3 of \cite{Tay3} (this is the only case when the lifts considered locally at $q$ do not have abelian image). Otherwise we choose a $\sigma=\sigma_q$ and a $\tau=\tau_q$ such that
$\rhobar(\sigma)$ is $$\left( \begin{array}{cc}
     r  & b  \\
    0 & r
\end{array} \right)$$ (note that by the relation $\sigma\tau\sigma^{-1}=\tau^q$, the characteristic polynomial 
of $\rhobar(\sigma)$ is forced to have double roots as $q$ is 1 mod $p$ and 
$\rhobar$ is tamely ramified at $q$), and $\rhobar(\tau)$ is the matrix $$\left( \begin{array}{cc}
     1  & a  \\
    0 & 1
\end{array} \right).$$  We want to construct infinitely many $O$-valued lifts $\rho_q$ of determinant 
$\vep\chi_p^{k(\rhob)-1}\eta_q^i|_{D_q}$, and such that
$\rho_q|_{I_q}$ is of the from
$$\left( \begin{array}{cc}
     \chi'  & *  \\
    0 & 1
\end{array} \right).$$ We seek $$\rho_q(\sigma)=\left( \begin{array}{cc}
     \alpha  & \gamma  \\
    0 & \beta
\end{array} \right),$$ say $A$, and $$\rho_q(\tau)=\left( \begin{array}{cc}
     \chi'(\tau)  & a'  \\
    0 & 1
\end{array} \right),$$ say $B$. 
Here  $\alpha\beta=\vep(\sigma)\chi_p^{k(\rhobar)-1}(\sigma)\eta_q^i(\sigma)$, $\alpha,\beta$ reduce to $r$, $a'$ reduces to $a$ (and hence is a unit), $\gamma$ 
to $b$. 
We explicitly produce these lifts by the following calculation. (The version 
of the calculation we present here is suggested by B\"ockle.)
We can assume by conjugation that $a$ and $a'$ are equal to one. 
As the order of $B$ divide $q-1$ (as $\chi'$ is not trivial of order dividing $q-1$), the relation 
$\rho_q(\sigma)\rho_q(\tau)\rho_q(\sigma)^{-1}=\rho_q(\tau)^q$
is therefore equivalent to $AB=BA$, which yields
     $\alpha-\beta=\gamma(\chi'(\tau)-1)$
  as the only relation. Combining this with the equation 
$\alpha\beta=\psi:=\epsilon(\sigma)
\chi_p^{k(\rhobar)-1}(\sigma)\eta_q^i(\sigma)$
  gives the quadratic equation
     $\beta^2-\beta\gamma(\chi'(\tau)-1)-\psi=0$ for $\beta$.
  Since $\chi'(\tau)-1$ lies in the maximal ideal, it follows easily that 
  for each $\gamma$ (reducing to $b$)
  there is a unique solution $\beta$ that is congruent to $r$ 
  modulo the maximal ideal of $ O$, and hence there is a unique $\alpha$ depending on $\gamma$.
  Thus one has a 1-parameter family (in $\gamma$) of lifts 
$\rho_q$ of the required type, thus proving that $R_q$ is smooth over $O$ of relative dimension 1.

\vspace{3mm}

After these 2 facts, from Section \ref{liftings}, we deduce that
the global deformation ring $R$ which parametrises (equivalence classes of) lifts of $\rhob$

- of the given determinant $\vep\chi_p^{k(\rhobar)-1}\eta_q^i$, 

- that at $q$ are of the given form,

- are minimal at primes $\neq q,p$ in the sense of the section above, 

- at $p$ the lift is crystalline of weight $k(\rhobar)$ 
(and thus finite flat when the residual representation has weight 2),

\noindent is a finite flat complete intersection (ffci) over $O$, and hence we get a lift of the desired kind.

\end{proof}

\noindent{\bf Remark:} We will apply this proposition only when $\rhobar$ has weight 2, and 
is unramified outside $p,q$, and the lifts that need to be constructed
have non-trivial nebentype at $q$. In the proposition, 
when $\rhobar$ is ordinary we could also have allowed the 
deformations at $p$ to be minimal of weight 2, and hence 
Barsotti-Tate over $\Q_p(\mu_p)$ when the weight is not $p+1$, and semistable of weight 2 otherwise.

As we have seen, the proofs of Proposition \ref{p} and \ref{q} follow easily from the method of proof of Theorem 2.2 of \cite{KW} after 
some computations of local deformation rings. 
The point is that because of the method of Section 2 of \cite{KW}, we can prove results about congruences of Galois representations to parallel many of the results known for congruences between modular forms.  The method should allow one to prove in many more cases the analog of the results in \cite{DT} and
\cite{inv}, these are level raising results for modular forms, for Galois representations: this is reduced 
to some local computation. The local computations in the ``$(p,p)$ case'' (needed for the analog of \cite{inv}) are likely to be involved, while those in the $\ell \neq p$ case may be easier in many cases, and could be deduced for instance from \cite{Boe} when locally the residual representation at the place $\ell$ is not scalar. 

\section{Compatible systems}\label{comp}

We explain how to make the lifts $\rho$ of Proposition \ref{p} and \ref{q}
part of a compatible system. As in Section 3 of
\cite{KW}, the proof uses the method of \cite{Tay2} (see proof of Theorem 6.6 of \cite{Tay2}), and the refinements in \cite{D} and \cite{Wint}.

\begin{prop}\label{c}
 
(i)  Assume $\rhobar$ is as in Proposition \ref{p}: so it is ordinary of weight $2 \leq k(\rhobar) \leq p+1$ and $k(\rhobar) \neq p$, and $p>3$.
Given a minimal lift $\rho$  of $\rhobar$ as in Proposition \ref{p}, there is a
(weakly) compatible system $(\rho_{\lambda})$ where $\lambda$ runs through all
places of a number field $E$, and $\rho$ is a member of the
compatible system above the prime of $p$ fixed by $\iota_p$. Further 
$\rho_{\lambda}$ for $\lambda$ above the prime $\ell$ ($>2$) fixed by $\iota_\ell$, that is not ramified in $\rhobar$ (and hence $\neq p$), 
the representation is Barsotti-Tate at $\ell$, unramified outside the primes ramified in $\rhobar$, and 
the inertial Weil-Deligne (WD) parameter at $p$ of $\rho_{\lambda}$ 
(i.e., if $(\tau,N)$ is the WD parameter, with $\tau$ a $F$-semisimple representation of the Weil group and $N$ 
a nilpotent matrix, we consider only $(\tau|_{I_q},N)$) is the same as that of
$\rho$ at $p$.
 
(ii) Now we assume $\rhob$ is as in Proposition \ref{q}: thus 
$\rhobar$ does not have solvable image (but $p=3$ is allowed), it has the behaviour at a prime $q$ as in  Proposition \ref{q}, but we make the additional assumption that 
$k(\rhobar)=2$. 
Consider a minimal lift $\rho$ as in Proposition \ref{q}, thus $\rho$
is Barsotti-Tate at $p$, and we assume that the nebentype at $q$, 
$\chi\eta_q^i$, is non-trivial. The inertial parameter at $q$ of $\rho$ is $(\omega_q^j \oplus 1,0)$ ($ 1 \leq j \leq q-2$) where 
$\omega_q^j:=\chi\eta_q^i$. There is a (weakly) compatible system $(\rho_{\lambda})$ where $\lambda$ runs
through all places of a number field $E$, and $\rho$ is a member of
the compatible system at the place above  $p$ fixed by $\iota_p$. Further $\rho_{\lambda}$ for
$\lambda$ a prime above $q$ that is determined by $\iota_q$, is unramified at all primes $\neq
p,q$ outside which $\rhobar$ is unramified, is unramified at $p$, and at $q$
is Barsotti-Tate over $\Q_q(\mu_q)$.  Assume that $\overline{\rho_\lambda}$
has non-solvable image. Then $\overline{\rho_\lambda}$ has weight $j+2$, 
or its twist by ${\overline \chi_q}^{-j}$
has weight $q+1-j$. 
\end{prop}

\begin{proof}
This follows from the results in \cite{Tay1} and \cite{Tay2}, using the arguments in Section 3 of
\cite{KW} (see Theorem 3.1 of \cite{KW}).  

The existence of a weakly compatible system $(\rho_{\lambda})$, of which $\rho$ is a member, follows easily from the proof of Theorem 6.6 of \cite{Tay2} (which uses Brauer's theorem on writing representations of finite groups as a virtual sum of representations induced from characters of solvable subgroups, and base change results of Arthur and Clozel in \cite{AC}). ({\small We recall the construction of Taylor in \cite{Tay1} and \cite{Tay2}. In both (i) and (ii) we may assume that the image of $\rhobar$ is non-solvable. Then the lift $\rho$ constructed is such that there is a totally real field $F$, Galois over $\Q$, such that $\rho|_{G_F}$ arises from a holomorphic, cuspidal automorphic representation $\pi$ of $\GL_2(\A_F)$ with respect to the embedding $\iota_p$. Using Brauer's theorem
we get subextensions $F_i$ of $F$ 
such that $G_i=\Galois (F/F_i)$ is solvable, characters $\chi_i$ of $G_i$ with values in $\overline \Q$ (that we embed in $\overline \Q_p$ using $\iota_p$), such that $1_{G}=\sum_{G_i}n_i{\rm Ind}_{G_i}^{G} \chi_i$. Using \cite{AC} we also
get holomorphic cuspidal automorphic representations $\pi_i$ of $\GL_2(\A_{F_i})$ such that
if $\rho_{\pi_i,\iota_p}$ is the representation of $G_{F_i}$ corresponding to $\pi_i$ w.r.t. $\iota_p$, then $\rho_{\pi_i,\iota_p}=\rho|_{G_{F_i}}$. Thus $\rho=\sum_{G_i}n_i{\rm Ind}_{G_{F_i}}^{G_\Q} \chi_i\otimes \rho_{\pi_i,\iota_p}$.
Now for any prime $\ell$ and any embedding 
$\iota:\overline \Q \rightarrow \overline \Q_\ell$, we define  the virtual representation
$\rho_{\iota}=\sum_{G_i}n_i{\rm Ind}_{G_{F_i}}^{G_\Q} \chi_i\otimes \rho_{\pi_i,\iota}$ of $\Gal$
with the $\chi_i$'s now regarded as $\ell$-adic characters via the embedding  $\iota$. We check that $\rho_{\iota}$ is a true representation by computing its inner product. The representations $\rho_\iota$ together constitute the weakly compatible system we seek.})

Thus we concentrate below on proving some of the finer ramification properties claimed for $\rho_{\lambda}$
at the prime of the same residue characteristic as $\lambda$.

For instance the property at $p$ asserted in (i) follows directly from
Taylor's results in \cite{Tay1} (see Lemma 1.5 and Corollary 1.7 of it) which show that $\rho|_{G_F}$ is modular for some $F$ that is unramified at $p$, i.e., it arises (w.r.t.
the embedding $\iota_p$) from a  Hilbert modular form $f$ for $F$.
In fact \cite{Tay1} also shows that $f$ is ordinary at all primes above $p$ (w.r.t. $\iota_p$), and at such primes 
is either principal series of conductor dividing $p$ with ``nebentype'' $\omega_p^{k-2}$ (i.e.,
the automorphic representation corresponding to $f$ is such that at primes $\wp$ above $p$ the corresponding
local component is the principal series $\pi(\psi_1,\psi_2)$, with $\psi_1$ restricted to the units given by the character  corresponding to  $\omega_p^{k-2}$ by local class field theory and $\psi_2$ unramified), or Steinberg at $p$ (the latter only in the case $k(\rhobar)=p+1$). The reader may also consult proof of Proposition 2.5 of \cite{KW} for more details, especially in the ordinary case.

We turn to proving (ii).
We know by \cite{Tay1} and \cite{Tay2} that there is a totally real field $F$
(which we may assume to be Galois over $\Q$) 
over which $\rho|_{G_F}$ is modular. 
If $F'\subset F$ is the fixed field of a  decomposition group above $q$, then there is a prime $Q$ above $q$ in $F'$ which is a split prime. We deduce  using \cite{AC}, \cite{Carayol} and \cite{Tay}, 
that $\rho|_{G_F'}$, and hence 
$\rho_\lambda|_{G_{F'}}$,  arises from a Hilbert modular form $f$ which locally at 
$Q$ is a ramified principal series of conductor $Q$, whose nebentypus at $Q$ is $\omega_q^j:=\chi\eta_q^i$. This also proves the assertion about the inertial parameter of $\rho$ at $q$.

We now focus on the properties of $\rho_{\lambda}$ asserted at $q$
(as outside $q$ the proof is the same as in Section 3 of \cite{KW}).(We may assume that $F$ is of even degree over $\Q$ as we are done otherwise by \cite{Saito}.) We first prove that $\rho_\lambda$ when restricted to
$\Q_q(\mu_q)$ is Barsotti-Tate (which uses the assumption
that $\chi\eta_q^i$ is non-trivial). To see this, we work over a totally real field 
$F''$ that is a solvable extension of $\Q$ which when completed at all 
places above $q$ is $\Q_q(\mu_q)$ (see Lemma 2.2 of \cite{Tay3}), 
and apply to $(\rho_\lambda|_{G_{F''}})$ the same arguments as those in proof of Theorem 3.1 of \cite{KW} for $(\rho_\lambda)$, but instead of using 
Th\'eor\`eme 1 (ii) of  \cite{Breuil-cong} we use 
a result of \cite{Raynaud} (see Proposition 2.3.1), which we may 
as we are in the weight 2 case. 
({\small Recall that $\rho|_{G_{F'}}$ arises from a weight 2 Hilbert modular form $f$ for $GL_2(\A_{F'})$ that when base changed to the composite $F'F''$ of $F'$ and $F''$, which is a solvable extension of $F'$, becomes
unramified at places above $q$. This uses results of \cite{Tay1}, \cite{Tay2}, 
\cite{Carayol}, \cite{Tay}, \cite{AC} as in Section 3 of \cite{KW}.
This then gives the required statement by deducing that 
$\rho_{\lambda}|_{G_{F''F'}}$ mod $\lambda^n$ is finite flat at primes above $q$ by Th\'eor\`eme 1 (i) of \cite{Breuil-cong},
and then using  Proposition 2.3.1 of \cite{Raynaud}, instead of using Th\'eor\`eme 1 (ii) of \cite{Breuil-cong} which got used in proof of Theorem 3.1 of \cite{KW}. Note that the completion at a prime above $Q$ of $F'F''$ is $\Q_q(\mu_q)$.}) 

Now assume that the image of $\overline \rho_\lambda$ is non-solvable.
It is not hard to see that 
${\overline \rho_\lambda}|_{G_{F'}}$, which we know is modular, also arises from  a Hilbert modular 
form $f'$, congruent to $f$ mod place fixed by $\iota_q$, that is square integrable at a finite place $\alpha$ and at $Q$ of conductor $Q$ and of nebentypus $\chi\eta_q^i$ at $Q$. 
({\small To get such a $f'$ congruent to $f$ is standard: 
We use the level raising techniques of \cite{Tay}, using the proof of Theorem 2 of {\it loc.\ cit.\ } to find a 
$\alpha$, prime to $q$,  such that using notation of Theorem 1 of {\it loc.\ cit.\ } the valuation under $\iota_p$ of $\Theta_f(T_{\alpha}^2-S_{\alpha}(\N \alpha +1)^2)$
is bigger than that of $E_f((\N \alpha +1))$ (for instance $\N$ stands for the norm from $F'$ to $\Q$). Then we prove the existence of the desired $f'$ which is square integrable at a finite place $\alpha$ and at $Q$ is fixed by $U_1(Q)$ which gives rise to ${\overline \rho_\lambda}|_{G_{F'}}$, using  the proof of Theorem 1 of {\it loc.\ cit.\ } with $\lambda=\alpha$ in the notation there.
That $f'$ has nebentype $\omega_q^j=\chi\eta_q^i$ at $Q$ follows by considering determinants and central characters.})

Then we use the results of \cite{Saito} to conclude that at $Q$ the inertial WD parameter of $\rho_\lambda|_{G_{F'}}$ 
is the same as that of $\rho|_{G_{F'}}$ which we know to be $(\omega_q^j \oplus 1, 0)$
(note that as $Q$ is a split prime of 
$F'$, local information at $Q$ of $\rho_\lambda|_{G_{F'}}$ gives the local information of $\rho_\lambda$
at $q$).  After this, to get the information about 
weights, we use Proposition 6.1.1 of \cite{BM} and Theorem 6.11 of \cite{Savitt1},
as we have that $\rho_\lambda|_{G_{F'}}$
is irreducible by Lemma 2.6 of \cite{KW}.
\end{proof}


\section{Chebyshev's estimates on primes}\label{cheb}

In the proof of Theorem \ref{main} below, 
we will need some estimates on prime numbers proven by Chebyshev (we learnt about 
the following precise form of Chebyshev's estimate from a message of J-P.~Serre, and from R.~Ramakrishna). 

If $\pi(x)$ is the prime counting function, then if $x >30$
$$A({ x \over {\rm log}(x)}) \leq \pi(x) \leq B({ x \over {\rm log}(x)})$$ where $A=0.921...$ and 
$B \over A$ is ${6 \over 5}=1.2$ (see \cite{Chebyshev1} and \cite{Chebyshev2}, and also page 21 of \cite{Ellison}). From this we easily deduce that
if we fix a real number $a>1.2$, and denote by $p_n$ the $n$th prime which we
assume $> {\rm max}(30,a^{{6} \over {5a-6}})$, then $p_{n+1} \leq ap_n$. 

In the arguments below this estimate will be relevant for particular values of $a$. Given $p_n>2$ we consider an odd (prime power) 
divisor $\ell^r=2m+1$ of 
$P_{n+1}-1$, where $P_{n+1}$ is either $p_{n+1}$, or $p_{n+2}$ if $p_{n+1}$ is a Fermat prime,  
and divide the integers in the interval $[0,P_{n+1}-1]$ into blocks of size ${{P_{n+1}-1} \over {\ell^r}}$. We need to  ensure that
$p_n+1 \geq {\rm max}({{m+1} \over {2m+1}}({P_{n+1}-1}  )+2,
P_{n+1}-({{m} \over {2m+1}}({P_{n+1}-1}  )))$. A computation shows that
this is ensured by requiring that
$${{P_{n+1}}\over p_n} \leq {2m+1 \over {m+1}} - ({m \over {m+1}})({1 \over {p_n}}).$$
We consider $p_n \geq 31$.

An inspection shows that there is always a $P_{n+1}$ as required up to $p_n=1000$: use $a$ in the Chebyshev estimate to be ${44 \over 30}={3 \over 2} -{1 \over 30}$ and rule out Fermat primes causing problems in that range (the only ones are $5,17,257$ and for $p_n=251$, the prime preceding $257$, we can use $P_{n+1}=263$).
After that the Chebyshev estimate used with $a=\sqrt{1.499}$ (note $1.499={3 \over 2} - {1 \over 1000}$, and that after $3,5$ no two successive primes can both be Fermat primes) 
gives the existence of $P_{n+1}$ of the required kind if $p_n$ is at least $21591$. For primes $1000<p_n \leq 21591$ we can again get the desired prime $P_{n+1}$ using $a={44 \over 30}$ as there are no Fermat primes in this range.
Thus we conclude that we always have a $P_{n+1}$ as desired 
once $p_n \geq  31$.

(We thank Bo-Hae Im and Faheem Mitha for help with the calculations of this section.)

\section{Four lemmas}

We have the following lemma which we owe, 
together with its proof which we reproduce {\it verbatim}, to Wintenberger
(it was in an early version of \cite{KW}). Below by a mod $p$ representation $\rhobar$ being semistable
we mean that for all the primes ramified in $\rhobar$ that are $\neq p$, the 
ramification is unipotent, and by $\rhobar$ being dihedral we mean its projective image is a dihedral group.

\begin{lemma}\label{dihedral}(Wintenberger)
 (i) A mod $p$ dihedral, semistable representation $\rhobar$
which is odd and irreducible, with $p$ odd, exists if  and only if $p$ is 3 mod $4$, and
the class number of the imaginary quadratic field $\Q(\sqrt{ -p })$
is non-trivial (i.e., $p$ is 3 mod 4 and $p \neq 3, 7, 11, 19,
43, 67,163$). Such a $\rhobar$ is split at $p$ and has a twist with Serre weight
$k={ {p+1} \over 2}$.
 (ii) If $\rhob$ is a dihedral representation  induced from the quadratic
subfield of $\Q(\mu_p)$,
and $\rhob$ is locally irreducible at $p$ of weight at most $p+1$, then its
weight is $p+3 \over 2$.
 
\end{lemma}

\begin{proof} (i) First note that such a dihedral $\rhobar$
is unramified outside $p$. Let the projective
image $H$ of $\rhobar$ be a dihedral group $D_{2t}$ of order $2t$, $t>1$.
As $H$ has order prime to $p$, the image of inertia $I_p$ in $H$ is cyclic
: denote by $\tau _p$ a generator of the image of $I_p$ in $H$, and by $i_p$ the order of $\tau_p$.
Note that $t$ is odd. For if $t$ were even, $D_{2t}$ would have
a quotient which is $\Z/2\Z \times \Z/2\Z$ which would give
a character of order $2$ unramified everywhere.
As $t$ is odd, $D_{2t}$ has only one
character of order $2$, say $\eta$. It corresponds to a quadratic
field $K$ which is ramified at $p$. As $\eta (\tau _p )=-1$,
and $t$ is odd, $\tau _p$ is of order $2$.
Let $L$ be the field cut up by the representation
of $\Gal$ in $H$. Then $L/K$ is unramified everywhere.

Let $c\in \Gal$ be a complex conjugation. As $\rhob$ is odd,
$c$ has non trivial image $\gamma$  in $H$. As $t$ is odd,
$\eta (\gamma)\neq 1$ and
the quadratic field  $K$ is imaginary. This implies
$p \equiv 3$ modulo $4$, and that $t$ is a divisor of $h_K$,
the class number of $K$ (which is known to be odd). Let $\Delta _p$ be the image of the decomposition group
at $p$ in $H$. As  $\Delta _p$ is a quotient of the Galois
group of $\Q_{p,unr}(\sqrt{p} )/\Q_p$ (wher $\Q_{p,nr}$ is the maximal unramified extension of $\Q_p$), $\Delta _p$ is
abelian. As the subgroup $\Gamma_p$ of $H$ of order $2$
generated by  $\tau _p$ is its own centralizer ($t$ is odd),
one sees that $\Delta _p=\Gamma_p$. This implies that $\rhobar(D_p)$ is abelian.Thus $I_p$ acts under
$\rhob$ via characters of level $1$.
Let us assume that we have twisted $\rhob$ by the suitable
power of the cyclotomic character so that
the Serre weight $k$ of $\rho$  satisfies $2\leq k\leq p+1$.
Then $I_p$ acts with characters $1$ and $\chi_p ^{k-1}$.
As the image of $I_p$ in $H$ is of order $2$, one has $k-1= {p-1 \over 2}$
and $k={p+1 \over 2}$.

(ii) This follows because the image of inertia $I_p$ in the projective 
image of $\rhob$, which we know to be cyclic, is forced to be of order 2.
\end{proof}

\noindent{\bf Remark:} We use part (i) of Lemma \ref{dihedral} only to
deduce that a dihedral $\rhobar$ which is unramified outside $p$ is up to twist
ordinary at $p$. Note that once we know that in the proof $t$ is odd, we can also see this by deducing that locally at $p$
the fixed field of the kernel of $\rhobar$ cannot have an unramified quadratic subfield.
\vspace{3mm}

We also have a useful lemma which follows immediately from Serre's definition of weights:

\begin{lemma}\label{weights}
  Assume $\tau:G_{\Q_p} \rightarrow GL_2(\F)$ is such that its weight $k 
\neq 2$ is $<p$. Then if $\tau$ is irreducible,
  and $k-1=p-k'$ with $k'$ non-negative, $\tau \otimes \overline{\chi_p}^{k'}$ has weight $k'+2=p+3-k$. If $\tau|_{I_p}$ is
  split, $\tau \otimes \overline{\chi_p}^{1-k}$ has weight $p+1-k$.
\end{lemma}

The following lemma is crucial for us and is due to Breuil and M\'ezard, \cite{BM}, Proposition 6.1.1, and Savitt, \cite{Savitt1}, Theorem 6.11:

\begin{lemma}\label{breuil}
If a representation $\rho: G_{\Q_p} \rightarrow GL_2({ O})$ 
becomes Barsotti-Tate over $\Q_p(\mu_p)$, then if residually the representation is reducible, then $\rho$ itself is reducible. More precisely, a twist of
$\rho$ by some power of $\omega_p$ (the Teichmuller lift of $\overline \chi_p$) is ordinary.
\end{lemma}

\begin{proof}
This follows from Proposition 6.1.1 of \cite{BM}, and Theorem 6.11 of \cite{Savitt1} as then we know that some twist of $\rho$ by a power of $\omega_p$ will be ordinary (i.e., in case (i) or (iii) of {\it loc.\ cit.\ } which note are related to each other by twists, and are Cartier dual).
\end{proof}

\vspace{3mm}

\noindent{\bf Remark:} The particular case of the lemma above for representations arising from modular forms is in Proposition 8.13 of \cite{Gross} or Section 6 of \cite{Edix}.
{\small A non-ordinary newform in $S_2(\Gamma_1(Np))$
which gives rise to a  $\rhobar$ that is reducible at $p$ is such that its $p$th Hecke eigenvalue $a_p$ has valuation 1 and its complex conjugate form is ordinary.} Further, if a mod $p$ newform in $S_2(\Gamma_1(Np))$ ($(N,p)=1$) has nebentype at $p$ given by 
$\overline{\chi_p}^j$ ($0 \leq j \leq p-2$), then the associated mod $p$ representation either has Serre weight $j+2$ or its twist by $\overline{\chi_p}^{-j}$ has weight $p+1-j$. (We are grateful to Bas Edixhoven for helpful correspondence about this.)

\vspace{3mm}

The following lemma is a simple consequence of modularity lifting theorems 
of \cite{Wiles}, \cite{TW}, \cite{SW2}, \cite{SW3}, \cite{DFG}, \cite{Tay1},\cite{Tay2} and the work in 
\cite{BLZ}.

\begin{lemma}\label{rubbish}
  Let $p$ be an odd prime. If $\rho:G_\Q \rightarrow GL_2({O})$ is an 
  irreducible $p$-adic representation unramified outside $p$ and at $p$ 
  crystalline of Hodge Tate weights $(k-1,0)$ with $k$ even and $2 \leq k \leq p+1$,
  and if the residual representation $\rhobar$ is modular,  
  then $\rho$ arises from  $S_k(SL_2(\Z))$.
\end{lemma}

\begin{proof}
  Using Lemma \ref{dihedral} we see that if $\rhobar$ is irreducible, and restricted to $G_{ \Q(\sqrt{{(-1)}^{p-1 \over 2}}p)  }$
  is reducible, then $\rhobar|_{I_p}$ is ordinary and distinguished.
  Thus the only case of the lemma not covered in the literature (see \cite{Wiles}, \cite{TW}, \cite{SW2}, \cite{SW3},\cite{DFG})
  is when $\rho$ is possibly of weight $p+1$ and is not ordinary. In fact this case does not occur, as using Corollary
  4.1.3 of \cite{BLZ}, we see that in this case the residual representation is 
such that $\rhobar|_{D_p}$ is irreducible and of Serre weight 2. But this contradicts the fact that
  there is no Serre-type $\rhobar$ of level 1 and weight 2 proved in Theorem 4.1 of \cite{KW}.
\end{proof}

\noindent{\bf Remark:} Using \cite{BLZ} it must be easy to prove a lifting theorem when at $p$ the lifts are crystalline of weight $p+1$ even without assuming that the lift is unramified outside $p$. As we do not need it in this paper, we do not do this. Note that oddness of the residual representation forces $k$ to be even (see Proposition 2 and Corollary 2 of \cite{Wint}).

\begin{cor}\label{trivial}
  (i) Given an odd  prime $p$, if all 2-dimensional, 
  mod $p$, odd, irreducible representations $\rhobar$
  of a given weight $k(\rhobar)=k \leq p+1$, and unramified
  outside $p$, are known to be modular,
  then for any prime $q \geq k-1$, all 2-dimensional, mod $q$, odd, irreducible representations $\rhobar'$
  of weight $k(\rhobar')=k$, and unramified outside $q$, are modular.
  
(ii) If the level 1 case of Serre's conjecture is known for a prime $p>2$,
  then for any prime $q$ the level 1 case of Serre's conjectures is known
  for all 2-dimensional, mod $q$, odd, irreducible representations $\rhobar$
  of weight $k(\rhobar) \leq p+1$.
\end{cor}

\begin{proof} This is proved implicitly in \cite{KW}.
It is enough to prove the first statement. We may assume $q>3$.
By Theorem 2.1 of \cite{KW} lift $\rhobar'$ to a representation $\rho'$ which is unramified outside $q$ and crystalline at $q$ of weight $k$. By Theorem 3.1
of {\em loc.\ cit.\ }, $\rho'$ is part of a compatible system $(\rho'_\lambda)$
such that at a place $\lambda$ above $p$ the representation is unramified outside $p$ and crystalline of weight $k$ at $p$. The hypothesis of the corollary and Lemma \ref{rubbish}, then gives the modularity of $(\rho'_\lambda)$ and hence that of $\rhobar'$.
\end{proof}

\section{Proof of Theorem \ref{main}}

\subsection{Small weights}

We now prove our main theorem up to weights 32. (We do not really have to do so many weights before giving the general argument, 
but this seems good preparation for that and also verification of the general strategy in concrete cases.) When we choose a place $\lambda$ above a prime $\ell$ this is always chosen to be with respect to the
embedding $\iota_\ell: \aQ \hookrightarrow \aQl$  fixed once and for all a while ago.

For weights 2,4,6 this has already been proved in \cite{KW}.  Although the cases of weights 8,12,14 are also done in \cite{KW}
we redo them to illustrate that after weight 6 our inductive method takes over.
Our arguments below prove that the $\rhobar$ being considered arises from some level and weight, and then that it arises from weight $k(\rhobar)$ and level $N(\rhobar)=1$ follows from the results in \cite{Ribet}, \cite{Edix}.
(Also note that the Serre weight in all cases considered is even and so 
the residual mod $p$ representations when restricted to $I_p$ are never 
scalar. We often use the fact implicitly below that if $\rhobar$ is modular so is any twist of it by an abelian character, and the same fact for twists of $p$-adic representations $\rho$ by finite order characters.)

-- Consider the
weight 8 case. It's enough to prove using Corollary \ref{trivial}
that an irreducible $2$-dimensional, mod $7$ representation $\rhob$ of level $1$ and
weight $8$ is modular. (If the image is solvable we are done.) Such a representation is ordinary at $7$. By Proposition 
\ref{p}, lift $\rhobar$ to a weight 2,
semistable at $7$, $7$-adic representation $\rho$ that is unramified outside $7$. Using part (i) of 
Proposition \ref{c} get a compatible system $(\rho_\lambda)$ and reduce it mod a prime above 3 determined by $\iota_3$ to get $\rhobar'$. Note that 
by Proposition \ref{c} (i), $k(\rhobar')=2$, and $\rhobar'$ is
unramified outside 3 and 7.  If $\rhobar'$ has solvable image we are done, as in that case we have a representation to which we can apply known 
modularity lifting results to conclude that the compatible system
$\rho_{\lambda}$ is modular (\cite{Wiles}, \cite{TW}, \cite{SW2}, \cite{SW3}) and hence so is $\rhobar$:
these modularity lifting results may be applied
as by part (ii)   of Lemma \ref{dihedral}, $\rhobar'$ cannot be both reducible when restricted
to $G_{\Q(\sqrt{-3})}$ and irreducible locally at $3$. Similarly 
if $\rhobar'$ is unramified at $7$ we are again done
as we know by page 710 of \cite{Serre2} that the residual mod $3$ representation is then reducible.
Otherwise use Proposition 
\ref{q}, to
get a 3-adic lift $\rho'$ of $\rhobar'$ with nebentype $\omega_7^2$ at 7 ($\omega_7^4$ would also work). Then use Proposition
\ref{c} to get a compatible system $(\rho'_{\lambda})$ with $\rho'$ the member of this compatible system at the place corresponding to $\iota_3$, 
and consider a residual representation $\rhob'_7$ arising from this system at a place 
$\lambda$ above 7 fixed by $\iota_7$. We know
by Lemma \ref{breuil} (Proposition 6.1.1 of \cite{BM})
that if the residual mod $7$ representation $\rhob'_7$ locally at $7$ is reducible then the $7$-adic
representation is also locally reducible. Thus in this case up to twisting 
by a power of $\omega_7$ we may
assume the $7$-adic representation is ordinary at $7$ by Proposition 6.1.1 of \cite{BM}: the representation will also be $I_7$-distinguished as the residual weights are even.

If $\rhob'_7$ has solvable image, and hence known to be modular, we are done by applying results of \cite{SW2}, \cite{SW3} and \cite{Kisin}, which can be applied as
in the ordinary cases the representation will be $I_7$-distinguished
(and that in the dihedral case the representation is ordinary at $7$ by part (i) of Lemma \ref{dihedral}), and we conclude that 
$(\rho'_{\lambda})$ is modular.

Now assume that $\rhob'_7$ has non-solvable image. We get a
residual representation whose Serre weight (up to twisting) is either 4
or $7+3-4=6$  by Proposition 6.1.1 of \cite{BM}, Theorems 6.11 and 6.12 of 
 \cite{Savitt1} (as explained in Proposition \ref{c} (ii)), 
and we know the residual modularity for such weights (in fact we also know that in these cases again that $\rhob'_7$ is reducible,
although for uniformity of treatment we do not use this).  
Now we again use modularity lifting results in \cite{SW2}, \cite{SW3}, \cite{Kisin} as we just used and conclude that 
$(\rho'_{\lambda})$ is modular.

Hence so is the residual representation 
$\rhobar'$, and hence by another application of modularity lifting theorem 
we conclude that the first compatible system $(\rho_{\lambda})$ is modular (as the compatible systems $(\rho_{\lambda})$ and $(\rho'_{\lambda})$ are linked at the place above 3 fixed by $\iota_3$), and hence so is $\rhobar$ (which in this case means that it does not exist!).

(Now we will be more succinct, and skimp some of the details before we get to the general arguments of the next section, as it's much of the same thing!)

-- Weights 10 and 12: Consider a representation $\rhob$ mod 11 that is irreducible and of
weight 10 or 12. Its enough by Corollary \ref{trivial} to prove this to be modular to conclude that for any prime at least 11, any $\rhob$ of level 1 and 
weight 10 or 12 is modular. (If the image is solvable we are done.) By Lemma \ref{weights} locally at 11 we may assume that the representation is (non semisimple and)
ordinary, and hence we get a weight 2 minimal lifting $\rho$ unramified outside 11. 
Then by Propositions \ref{p} and  \ref{c} make $\rho$ part of a 
compatible system $(\rho_{\lambda})$, 
such that the representation in the system corresponding to $\iota_{11}$ is $\rho$,
and then consider the residual representation $\overline{\rho_{5}}$ 
at a prime $\lambda$ above 5 fixed by $\iota_5$.  We are done if the image
of $\overline{\rho_{5}}$ is solvable (or unramified at 11 and hence solvable by the weight 2, level 1 case proved in \cite{KW}) arguing just as in the previous case. 
Otherwise using Proposition \ref{q} construct lift of $\overline{\rho_{5}}$ with nebentype $\omega_{11}^4$ ($\omega_{11}^6$ would also work) at 11, and we can do this for either of the 2 weights 10 or 12 being considered.
Argue as before and get a compatible system $(\rho'_{\lambda})$ which residually 
at the  prime above 11 fixed by $\iota_{11}$ 
will have weight (up to twisting by some power of $\overline \chi_{11}$)
either 6 or 8 (the same weights if we had chosen lifting with nebentype $\omega_{11}^6$  when applying Proposition \ref{q}), and then we are done as before.

--Weights $14,16,18,20$: It will be enough to show that a mod $19$ representation $\rhob$ 
(irreducible, odd, $2$-dimensional of level one as always) of any of these weights is modular. (If the image is solvable we are done.) As we have 
dealt with weights up to $12$, we may assume by Lemma \ref{weights} that $\rhob$ is ordinary at $19$.
We again construct a weight 2 minimal lift $\rho$ of $\rhobar$ and get a compatible system $(\rho_{\lambda})$ (by Propositions \ref{p} and \ref{c}(i)) and consider a residual representation at the place above $3$
determined by $\iota_3$ (again we are done if the residual mod 3 representation has either solvable image or is unramified at 19), and apply Proposition \ref{q}
to get a lifting $\rho'$ of the mod 3 representation  that is unramified outside $3,19$, Barsotti-Tate at 3, and has nebentype
$\omega_{19}^8$ at 19 ($\omega_{19}^{10}$ would also work). Then make it part of a compatible system $(\rho'_{\lambda})$ (using Proposition \ref{c} (ii)) 
which residually at a prime above 19 will 
have weight (up to twisting by some power of $\overline \chi_{19}$) either 10 or 12 (the same if the nebentype $\omega_{19}^{10}$ had been chosen). Again we can argue as 
before, knowing Serre's conjecture in level 1 for weights 10 and 12,  and we are done.

-- Weights $22,24,26,28,30$: It will be enough to show that a mod 29 representation $\rhob$ (irreducible, odd, 2-dimensional of level one as always) of any of these weights is modular. (If the image is solvable we are done.) As we have dealt with weights up to 20, 
we may assume by Lemma \ref{weights} that $\rhob$ is ordinary at 29.
This time we lift $\rhobar$ to a minimal lift $\rho$ of weight 2
(by Proposition \ref{p}), make $\rho$ part of 
a  compatible system $(\rho_{\lambda})$ 
(by Proposition \ref{c}(i)) and 
consider a residual representation at a prime above 7 determined by $\iota_7$
(that we may assume is ramified at $29$ and has non-solvable image). Apply Proposition \ref{q}
to get a lifting $\rho'$ of the mod $7$ representation  that is unramified outside $3,19$, Barsotti-Tate at $7$, and has nebentype
$\omega_{29}^{16}$ at 29 when the weight is one of 22,26,30, or nebentype
$\omega_{29}^{14}$ at 29 when the weight is 24 or 28. Using Proposition \ref{c} (ii)
make it part of a compatible system $(\rho'_{\lambda})$ which residually at a prime above 29 will have weight (up to twisting by some power of $\overline \chi_{29}$) 18 or 14 (if the 
weight of $\rhobar$ is one of  22,26,30), or weight 16 (if the weight of $\rhobar$ was either 24 or 28). Again we can argue as 
before, knowing Serre's conjecture in level 1 for weights 14,16,18,  and we are done.

-- Weight 32: It will be enough to show that a mod 31 representation $\rhob$ (irreducible, odd, 2-dimensional of level one as always) of weight 32 is modular. (If the image is solvable we are done.) It's the same argument as before. We use as a foil the prime 5, and in the end construct a weight 2 compatible system $(\rho'_{\lambda})$, whose modularity yields that of $\rhobar$, 
of nebentype $\omega_{31}^{16}$ at 31
such that at a prime above 31 the residual representation has weight (up to twisting by some power of $\overline \chi_{31}$)
either 18 or 16. As we know Serre's conjecture in level 1 for these weights we can conclude.

\subsection{The general argument}

Now we give the general argument (which the reader 
must surely have guessed the gist of). This, together with the previous section, will prove Theorem \ref{main}.

Assume we have proven the level 1 case of Serre's 
conjecture mod $p_n$ where $p_n$ is a prime $\geq 31$. This implies by Corollary \ref{trivial} (ii)
that for any prime $q$ we know Serre's conjecture for any level 1 
mod $q$  representation $\rhobar$ that is odd, irreducible, 2-dimensional
of weight $ k(\rhobar) \leq p_n+1$.
 
By the arguments in Section \ref{cheb}, we can find a prime $P_{n+1}>p_n$, 
which is not a Fermat prime, and an odd prime power $\ell^r=2m+1$ that divides $P_{n+1}-1$ exactly such that
\begin{equation}\label{*}
{{P_{n+1}} \over {p_n}} \leq 
{2m+1 \over m+1} - ({m \over m+1})({1 \over p_n}).
\end{equation} (This inequality holds for any integer $m \geq 1$.)
Consider any weight $k$ such that $p_n+2 \leq k \leq P_{n+1}+1$: 
given such a $k$ 
we would like
to prove that any 2-dimensional irreducible, odd, mod $P_{n+1}$ representation $\rhobar$ of $\Gal$ of level 1 and weight $k(\rhobar)=k$ is modular. 
By Corollary \ref{trivial}, this will prove Serre's level 1 conjecture
for any 2-dimensional, odd irreducible, mod $q$ representation 
of Serre weight $k$, where $q  \geq k-1$, and also the level 1 conjecture mod all primes $\leq P_{n+1}$ 
once we have done this for all weights $p_n+2 \leq k \leq P_{n+1}+1$. Then we continue with $P_{n+1}$, treating it like $p_n$ etc.
We denote $P_{n+1}$ by $p$ for notational 
simplicity. (This is the inductive method 
to attack the level 1 case of Serre's conjecture 
proposed in Theorem 5.1 of \cite{KW}.)

By  Lemma \ref{weights}, and as we know the conjecture for weights $\leq p_n+1$, and using (\ref{*}),
we can assume that $\rhobar$ is ordinary at $p$. 
Construct a minimal weight 2 
lift $\rho$ of $\rhobar$ using Proposition \ref{p}. (If $\rhobar$ is solvable we already know its modularity, so we may assume it is not, although we need not.) 
This lift is unramified outside $p$ and minimal of weight 2 at $p$.
Using part (i) of Proposition \ref{c} make it part of a compatible system $(\rho_{\lambda})$, so that for a place above $p$ fixed by the embedding $\iota_p$ the corresponding representation is $\rho$. 
Now reduce the system modulo a prime $\lambda$ above 
$\ell$ with $\lambda$ determined by the embedding $\iota_{\ell}$. 
Denote the residual  mod $\ell$ representaion by $\rhobar_\ell$, and note that by part (i) of Proposition \ref{c}, $k(\rhobar_\ell)=2$.
If $\rhobar_\ell$ had solvable image we would be done. This is because by Lemma \ref{dihedral} (i), 
in the case when the representation $\rhobar_\ell$ restricted to $G_{\Q(\sqrt{{(-1)}^{\ell-1 \over 2}\ell})}$ is reducible, 
$\rhobar_\ell|_{D_\ell}$ is ordinary (as $\ell>2$ and hence $\ell+3 \over 2$$>2$), 
and in the cases when either this happens or $\rhobar_\ell$ is reducible, and thus ordinary at $\ell$, 
$\rhobar_\ell$ restricted to $I_{\ell}$ is distinguished.  As $\ell$ is an odd prime, 
this allows us to apply the modularity lifting theorems 
proved in \cite{Wiles}, \cite{TW}, \cite{SW2} and \cite{SW3}, to prove the modularity of 
$\rho_\lambda$ and hence of $\rhobar$.
If the mod $\ell$ representation were unramified at 
$p$, using the result in \cite{KW} for weight 2, we 
would see that the representation is reducible and 
again we would be done.

Otherwise apply Proposition \ref{q}
(with $p$ of that proposition, the present $\ell$, and the $q$ there the present $p$!) 
to choose a lifting $\rho'$ of this residual mod $\ell$ representation $\rhobar_\ell$, unramified outside $\ell,p$, 
Barsotti-Tate at $\ell$,   and such  
that the nebentype at $p$ is $\omega_p^j$  
with $j$ in the interval $({m \over 2m+1} ({p-1}),{m+1 \over 2m+1} ({p-1})]$
(recall that we have set $\ell^r=2m+1$).
Using part (ii) of Proposition \ref{c}, make $\rho'$ part of a compatible system $(\rho'_{\lambda})$ and consider a place above $p$
determined by the embedding $\iota_p$. We get a $p$-adic representation $\rho''$ at that place 
that is unramified outside $p$, and such that
at $p$ the representation is Barsotti-Tate over (the degree $p-1 \over 2$ subfield of) $\Q_p(\mu_p)$ by part (ii) of Proposition \ref{c}. 

Now note that if the representation residually, say $\rhobar''$, were reducible locally at $p$, then $\rho''$ at $p$ itself would be ordinary up to twisting by 
a power of $\omega_p$ by 
Lemma \ref{breuil} above (which is a consequence of 
Proposition 6.1.1 of \cite{BM}, and Theorems 6.11  of 
 \cite{Savitt1}), and would also
have distinct characters on the diagonals when restricted to $I_p$ (as the weight is even), i.e., is $I_p$-distinguished. 

Further if $\rhobar''$ had solvable image, then  as by 
 Lemma \ref{dihedral} above we know that if $\rhobar''$ is reducible when restricted to $G_{\Q(\sqrt{{(-1)}^{p-1 \over 2}p})}$
then it's ordinary at $p$, and also $I_p$-distinguished, we can apply the modularity lifting theorems of 
\cite{SW2},\cite{SW3} and \cite{Kisin} to conclude the modularity of $\rho''$.
The main theorem of \cite{Kisin}  
is used in the non-ordinary case,
when we invoke a modularity lifting result when the lift is potentially Barsotti-Tate at $p$ over a tamely ramified extension, in fact over $\Q_p(\mu_p)$, and the residual representation is irreducible when restricted to
$G_{\Q(\sqrt{{(-1)}^{p-1 \over 2}}p)}$.

So we can now assume that $\rhobar''$ has non-solvable image. 
We see by the last line of part (ii)  of Proposition \ref{c},
that $k(\rhobar'')$  is either $j+2$ or $k(\rhobar'' \otimes \overline \chi_p^{-j})=p+1-j$. Note that by our choice of $j$, $j+2$ is contained 
in in the interval $({ m \over 2m+1} ({p-1})+2,({ m+1 \over 2m+1}) ({p-1})+2]$ and $p+1-j$
is in the interval
$[p+1-(({ m+1 \over 2m+1}) ({p-1})), p+1-(({ m \over 2m+1}) ({p-1})))$.
 But (\ref{*}) gives that both these intervals are contained in the interval $[2,p_n+1]$  as noted in Section \ref{cheb}.
Hence we know the modularity of $\rhobar''$ as we know by  hypothesis the modularity of irreducible, mod $p$, odd, 2-dimensional, unramified outside $p$ 
representations of weights $\leq p_n+1$.
We can now apply the modularity lifting theorems of 
\cite{Kisin} to conclude that $\rho''$ is modular.

Thus the compatible system $(\rho_\lambda')$ is modular. Note that $(\rho_\lambda)$ 
and $(\rho_\lambda')$ are linked at the prime above $\ell$ determined by $\iota_\ell$, i.e., 
at the prime 
above $\ell$ fixed by $\iota_\ell$, the residual representations 
arising from the 2 systems are isomorphic (and have non-solvable image). Thus applying the modularity lifting theorems which were the first ones to be proven, i.e., in \cite{Wiles} and \cite{TW}, 
we conclude that $(\rho_{\lambda})$ is modular. Hence $\rhobar$ is modular of some weight $k$ and level, and then of weight $k(\rhobar)$ and level 1 
by \cite{Ribet}, \cite{Edix}.

\vspace{3mm}

\noindent{\bf Remark:} It will be of interest to see, both in \cite{KW} and the present paper, if the residually solvable cases can also be handled internally by the method of the papers themselves, rather than using the fact that Serre's conjecture is known for them.
%

We also remark that the proof does use an auxiliary  prime $\ell \neq p$, but no ramification is introduced at this prime, unlike in the 
more traditional use of auxiliary primes. 

\section{Corollaries}\label{cors}

\subsection{Proof of corollaries}

The case of $q=2$ of Corollary \ref{cond} is dealt with in Theorem 4.1(ii) 
of \cite{KW}. (The case $p=3$ is excluded there, but is dealt with
by the remarks about $p=3$ in Section \ref{liftings}.) Note that in the case $q=2$ we are using results of the type proven in \cite{Fontaine}, \cite{BK}, \cite{Schoof1}.

Corollary \ref{cond} when $q \neq 2$ follows from Theorem \ref{main} on using the method of ``killing ramification'' in Section 5.2 of \cite{KW}. We give very briefly the proof.
First apply the minimal lifting result of \cite{KW}, Theorem 2.1 of Section 2,
to construct first a minimal $p$-adic lift $\rho$ of $\rhobar$ and then apply Theorem 3.1 of 
{\it loc.\ cit.\ } to get a compatible system $(\rho_\lambda)$
of which $\rho$ is a part.
Consider the prime above $q$ determined by $\iota_q$ and reduce the representation which is a member of $(\rho_\lambda)$ 
at this prime (this is the method of ``killing ramification'' of Section 5.2 of 
\cite{KW}). Residually we get a mod $q$ representation that is unramified outside $q$, for which we know Serre's conjecture by the main theorem of this paper. Then the lifting theorems recalled above, which may be applied because of  Lemma \ref{dihedral} and Lemma \ref{breuil} (the latter due to Breuil, M\'ezard and Savitt), imply that the compatible system $(\rho_\lambda)$ is modular and hence so is $\rhobar$.

Corollary \ref{finite} follows as we know that a mod $p$ irreducible 2-dimensional representation
of $G_\Q$ that arises from $S_k(SL_2(\Z))$ (for any integer $k\geq     2$) also arises from $S_2(\Gamma_1(p^2))$.

\subsection{Quantitative refinements?} 

Corollary \ref{finite} is a folk-lore consequence of Serre's conjecture and 
is there implicitly at the end of Tate's  article \cite{Tate}.
It will be of interest to get quantitative refinements of Corollary \ref{finite}. 

 After the corollary it is easy to see that for a fixed prime $p$ the number
$N(2,p)$ of
isomorphism classes of continuous semisimple
 odd representations $\rhob:G_\Q \rightarrow GL_2(\aFp)$
 that are unramified outside $p$ is bounded by $Cp^3$ for a constant $C$ independent of $p$. This is seen as by the main theorem of the paper we are counting the number of distinct level 1 mod $p$ Hecke eigensystems, and then using formulas for the dimension of $S_k(SL_2(\Z))$, and the fact that, up to twist by powers of $\overline \chi_p$, all mod $p$ forms have weight $\leq p+1$, we are done. 

We would guess that $N(2,p)$ is asymptotic to ${1 \over {48}} p^3$ with $p$.

Serre has pointed out that another easy corollary of Theorem \ref{main} is that any  $\rhob$  in Corollary \ref{finite} can be
written over a finite field  $\F_{p^e}$  with  $e \leq {\rm sup}(1,{p+1 \over 12})$.

\section{Some remarks}

The arguments in this paper were arrived at when thinking of how to extend the results in \cite{KW} 
to prove Serre's conjecture in the level 1 case for 
(finitely many) more weights, upon Schoof telling us that he could prove that all semistable abelian varieties over $\Q$ with good reduction outside 17 (resp., 19) are isogenous to powers of $J_0(17)$ (resp., $J_0(19)$). Using the method of \cite{KW} this immediately proved modularity of $\rhobar$ of level 1 and weights 18 and 20. With further thought, we could do weights 10 and 16. The trick we came up with to do weight 16
led to the method presented in this paper.

The trick to do weight 10, which we have not used in the
present paper, was as follows: it's enough to prove that a mod $11$ representation of weight 10 is modular, and enough to assume that it's ordinary at 11. Now we can get a minimal lifting of $\rhobar$ that is crystalline, and ordinary, of weight 20 at $11$ by using the arguments of proof of Theorem 2.1 of \cite{KW}, and then putting it in a compatible system and 
reducing it mod $19$ we deduce modularity by the weights $2$ and $20$ results.
The trick for weight 16 is subsumed, and superseded (as we at first still made
use of Schoof's result for the prime 19), in the 
method of the paper. 

The weaning away (after weight 6) 
from results classifying semistable abelian varieties 
$A$ with good reduction outside a prime $p$ was necessary, as such results can only be proven for small $p$. For large $p$ such $A$ may not be isogenous to products of $\GL_2$-type abelian varieties, and even those of the latter type cannot be proved to be factors of $J_0(p)$  by arguments that use discriminant bounds. For a related reason the method of Tate cannot be used
to prove Theorem \ref{main} for large primes $p$ (perhaps $p=5$ is the limit
of this method, even assuming the GRH: see \cite{Tate} and the discussion in \cite{Brue}).

In \cite{KW} a possible, and even plausible, path to the general case of Serre's conjecture was mapped out
in its last section. The method was to use induction on primes in 2 different ways. For the level 1 case (see Theorem 5.1 there), the induction was on the prime which was the residue characteristic, the starting point being the cases of the conjecture proved for $p=2,3$ by Tate and Serre. 
We have used this method in the present paper.
For the reduction to the level 1 case (see Theorem 5.2 there; this is the method of ``killing ramification''),
the induction was on the number of primes ramified in the residual representation, the starting point being the level 1 case.  

But the path seemed to be blocked by 
formidable obstacles. This paper reaches 
the level 1 case by sidestepping (and not overcoming any of) these technical obstacles. The method here, combined 
with the proposed method of reduction of the general case 
to the level 1 case in the last section of \cite{KW}, will probably be useful for general
level $N$, and ease some of the technical difficulties. 




\vspace{5mm}

\noindent{\bf Acknowledgements:} We are grateful to Ren\'e Schoof for 
telling us about his unpublished results 
for semistable abelian varieties over $\Q$ with good reduction outside 17 and 19: although in this paper we have 
not directly made use of these, our attempt to apply them to Serre's conjecture led to the method of the paper.
We would like to thank Gebhard B\"ockle for very helpful correspondence about Section 2 of the paper. 
We are grateful to Gebhard B\"ockle, Kevin Buzzard, Jean-Pierre Serre and Jean-Pierre Wintenberger for very helpful comments on the mansucript. We thank Christophe Breuil for pointing out the reference \cite{Savitt1}. We thank Jean-Pierre Serre for telling us about the history of his conjecture.  The author was partially supported
by NSF Grant DMS - 0355528.

\nocite{*}
\bibliographystyle{plain}

\end{document}